\documentclass[11pt]{article}
\usepackage{amsmath, amsthm}
\usepackage{latexsym}
\usepackage{graphicx}
\makeatletter
\newfont{\footsc}{cmcsc10 at 8truept}
\newfont{\footbf}{cmbx10 at 8truept}
\newfont{\footrm}{cmr10 at 10truept}
\makeatother
\newtheorem{theorem}{\bf Theorem}

\newtheorem{definition}{\bf Definition}

\begin{document}
\title{Stochastic Ordering of Exponential Family Distributions and Their Mixtures}

\author{Yaming Yu\\
\small Department of Statistics\\[-0.8ex]
\small University of California\\[-0.8ex] 
\small Irvine, CA 92697, USA\\[-0.8ex]
\small \texttt{yamingy@uci.edu}}

\date{}
\maketitle

\begin{abstract}
We investigate stochastic comparisons between exponential family distributions and their mixtures with 
respect to the usual stochastic order, the hazard rate order, the reversed hazard rate order, and the 
likelihood ratio order.  A general theorem based on the notion of {\it relative log-concavity} is shown 
to unify various specific results for the Poisson, binomial, negative binomial, and gamma distributions 
in recent literature.  By expressing a convolution of gamma distributions with arbitrary scale and shape 
parameters as a scale mixture of gamma distributions, we obtain comparison theorems concerning such 
convolutions that generalize some known results.  Analogous results on convolutions of negative 
binomial distributions are also discussed. 

\end{abstract}


\section{Stochastic orders and some general observations}
The study of stochastic orders has received attention in diverse areas including economics, operations research, reliability, 
and statistics (e.g., survival analysis).  For book-length treatments of both theory and applications, see Shaked and 
Shanthikumar (1994, 2007).  This paper is mainly concerned with four orders, namely the usual stochastic order $\leq_{st}$, 
the hazard rate order $\leq_{hr}$, the reversed hazard rate order $\leq_{rh}$, and the likelihood ratio order $\leq_{lr}$.  We 
recall the familiar definitions. 

\begin{definition}
Let $X$ and $Y$ be continuous random variables on $\mathbf{R}$ with probability density functions (pdfs), or discrete random 
variables on $\mathbf{Z}$ with probability mass functions (pmfs), $f(x)$ and $g(x)$, respectively.  Denote their respective 
cumulative distribution functions (cdfs) by $F(x)$ and $G(x)$.
\begin{itemize}
\item
$X$ is said to be smaller than $Y$ in the usual {\it stochastic order}, or $X\leq_{st} Y$, if $\bar{F}(x)\leq \bar{G}(x)$ for 
all $x$, where $\bar{F}(x)=1-F(x)$ and $\bar{G}(x)=1-G(x)$.
\item
$X$ is said to be smaller than $Y$ in the {\it hazard rate order}, or $X\leq_{hr} Y$, if $f(x)/\bar{F}(x)\geq 
g(x)/\bar{G}(x)$ for all $x$. 
\item
$X$ is said to be smaller than $Y$ in the {\it reversed hazard rate order}, or $X\leq_{rh} Y$, if $f(x)/F(x)\leq 
g(x)/G(x)$ for all $x$. 
\item
$X$ is said to be smaller than $Y$ in the {\it likelihood ratio order}, or $X\leq_{lr} Y$, if the likelihood ratio $f(x)/g(x)$ 
is a monotone decreasing function on the set $\{x:\ f(x)>0\ {\rm or}\ g(x)>0\}$. By convention $a/0=\infty$ whenever $a>0$.
\end{itemize}
\end{definition}
As is well-known, $X\leq_{lr} Y$ implies $X\leq_{hr} Y$ and $X\leq_{rh} Y$, either of 
which in turn implies $X\leq_{st} Y$.  Further basic properties of these orders can be found in Shaked and Shanthikumar (1994).

Despite their importance, to verify the relations $\leq_{st},\ \leq_{hr},\ \leq_{rh}$ or 
$\leq_{lr}$ can be nontrivial, e.g., when the relevant distributions are not in closed form.  This work provides some 
simple conditions that unify and generalize many results for specific distributions in recent literature.  The following {\it 
relative log-concavity order}, introduced by Whitt (1985) (see also Yu 2009), plays a critical role in the development. 

\begin{definition}
Let $X$ and $Y$ be continuous (discrete) random variables with pdfs (pmfs) $f(x)$ and $g(x)$ respectively.  We say $X$ is {\it 
log-concave relative to} $Y$, denoted $X\leq_{lc} Y$, if 
\begin{enumerate}
\item
the support of $X$, $supp(X)=\{x:\ f(x)>0\}$ and the support of $Y$, $supp(Y)=\{x:\ g(x)>0\}$ are both intervals on 
$\mathbf{R}$ ($\mathbf{Z}$);
\item
$supp(X)\subset supp(Y)$; 
\item
$\log(f(x)/g(x))$ is a concave function on $supp(X)$.
\end{enumerate} 
\end{definition}

The order $\leq_{lc}$ provides a way of deriving conditions that imply the four ground-level orders $\leq_{st},\ \leq_{hr},\ 
\leq_{rh}$, and $\leq_{lr}$.  This is analogous to gaining understanding of the monotonicity properties of a function by 
studying its second derivative.  We summarize some general observations below.  

\begin{theorem}
\label{gen}
Let random variables $X$ and $Y$ have pdfs $f(x)$ and $g(x)$ respectively, both supported on $(0, \infty)$.  Assume the log 
density ratio $l(x)=\log(f(x)/g(x))$ is continuous and moreover concave, i.e., $X\leq_{lc} Y$.  Then 
\begin{enumerate}
\item
$X\leq_{st} Y$ and $X\leq_{hr} Y$ are equivalent, and each holds if and only if $\lim_{x\downarrow 0} l(x)\geq 0$;
\item
assuming $l(x)$ is continuously differentiable, then $X\leq_{lr} Y$ and $X\leq_{rh} Y$ are equivalent, and each holds 
if and only if $\lim_{x\downarrow 0} l'(x)\leq 0$.
\end{enumerate}
\end{theorem}
\begin{proof}  Part 1).  Let $A=\{x:\ l(x)\geq 0\}=\{x:\ f(x)\geq g(x),\ x>0\}$.  Because $l(x)$ is concave, $A$ is an 
interval.  We first show that $X\leq_{st} Y$ is equivalent to $\lim_{x\downarrow 0} l(x)\geq 0$.  If $\lim_{x\downarrow 0} 
l(x)\geq 0$ then it is easy to see that the left end point of $A$ is 0.  That is, $f(x)-g(x)$ changes sign at most once from 
$+$ to $-$ as $x$ increases from $0$ to $\infty$; it follows that $F(x)-G(x)=\int_0^x (f(u)-g(u))\, {\rm d}u$ does not change 
sign at all, i.e., $F(x)\geq G(x)$ for all $x$, and $X\leq_{st} Y$ by definition.  Conversely, if $X\leq_{st} Y$ then 
$\int_0^x (f(u)-g(u))\, {\rm d}u\geq 0$ for all $x$, forcing the left end point of $A$ to be zero, which implies 
$\lim_{x\downarrow 0} l(x)\geq 0$.  Note that this limit exists by the concavity of $l(x)$. 

Concerning the hazard rate order, we only need to show $X\leq_{st} Y \Rightarrow X\leq_{hr} Y$ since the implication $X\leq_{hr} Y \Rightarrow X\leq_{st} Y$ is well-known.  By definition, if $X\leq_{st} Y$ then $\bar{F}(x)\leq \bar{G}(x)$ for all $x$.  Given $x_0>0$, if $f(x_0)\geq g(x_0)$, then $f(x_0)/\bar{F}(x_0)\geq g(x_0)/\bar{G}(x_0)$.  Otherwise $f(x_0)< g(x_0)$, i.e., $x_0\notin A$.  As before, since $X\leq_{st} Y$, the left end point of $A$ must be zero.  Hence $l(x)< 0$ for all $x\geq x_0$.  If there 
exist some $x_2>x_1\geq x_0$ such that $l(x_2)>l(x_1)$, then by the concavity of $l(x)$, for all $x\leq  x_1$ we have
$$l(x)\leq l(x_1)+(x-x_1)\frac{l(x_2)-l(x_1)}{x_2-x_1}< 0,$$
i.e., $l(x)<0$ for all $x$, a contradiction.  Thus $l(x)$ (or $f(x)/g(x)$) decreases on $[x_0, \infty)$, and consequently, 
\begin{align*}
\frac{f(x_0)}{\bar{F}(x_0)}&=\frac{f(x_0)}{\int_{x_0}^\infty f(u)\, {\rm d}u}\\
                           &\geq \frac{f(x_0)}{\int_{x_0}^\infty g(u)f(x_0)/g(x_0)\, {\rm d}u}\\
                           &=\frac{g(x_0)}{\bar{G}(x_0)}.
\end{align*}
That is, the hazard rate of $X$ is always greater than or equal to that of $Y$.

Part 2).  Note that $l'(x)$ decreases in $x$ since $l(x)$ is concave; therefore to ensure monotone density 
ratio, or $l'(x)\leq 0$ for all $x$, we only need $\lim_{x\downarrow 0} l'(x)\leq 0$.
That is, 
$$X\leq_{lr} Y\Longleftrightarrow \lim_{x\downarrow 0} l'(x)\leq 0.$$

Concerning the reversed hazard rate order, we only need to show $X\leq_{rh} Y\Rightarrow X\leq_{lr} Y$, since 
the implication $X\leq_{lr} Y\Rightarrow X\leq_{rh} Y$ is known.  Assume the contrary, i.e., $X\leq_{rh} Y$ but 
$X\not\leq_{lr} Y$.  Then, by the discussion above, $\lim_{x\downarrow 0} l'(x)>0$, and by continuity 
there exists $\epsilon>0$ such that $l'(x)>0$ for all $x\in (0, \epsilon]$.  That is, $f(x)/g(x)$ strictly increases on $x\in 
(0, \epsilon]$.  Thus 
\begin{align*}
\frac{f(\epsilon)}{F(\epsilon)} &=\frac{f(\epsilon)}{\int_0^\epsilon f(u)\, {\rm d}u}\\
                                &> \frac{f(\epsilon)}{\int_0^\epsilon g(u) f(\epsilon)/g(\epsilon)\, {\rm d}u}\\
                                &=\frac{g(\epsilon)}{G(\epsilon)}
\end{align*}
which contradicts the definition of $X\leq_{rh} Y$.
\end{proof}

A discrete version of Theorem \ref{gen} is
\begin{theorem}
\label{discrete}
Let random variables $X$ and $Y$ have pmfs $f(x)$ and $g(x)$ respectively, both supported on the same $\mathbf{Z}_+=\{0, 1, 
\ldots\}$ (or $\{0, 1,\ldots, n\}$ for some $n>0$).  Assume $X\leq_{lc} Y$.  Then 
\begin{enumerate}
\item
$X\leq_{st} Y$ and $X\leq_{hr} Y$ are equivalent, and each holds if and only if $f(0)/g(0)\geq 1$;
\item
$X\leq_{lr} Y$ and $X\leq_{rh} Y$ are equivalent, and each holds if and only if $f(1)/g(1)\leq f(0)/g(0)$.
\end{enumerate}
\end{theorem}

Basically, if $X\leq_{lc} Y$, then all of $X\leq_{st} Y,\ X\leq_{hr} Y,\ X\leq_{rh} Y$ and $X\leq_{lr} Y$ are determined by 
the behavior of $\Pr(X=x)/\Pr(Y=x)$ near the left end point $x=0$.  

{\bf Example 0.}  Let $Y\sim {\rm Bin}(n, p),\ p\in (0, 1)$, and $X=\sum_{i=1}^n B_i$, where $B_i$ are independent Bernoulli 
random variables, i.e., $\Pr(B_i=1)=1-\Pr(B_i=0)=p_i,\ i=1, \ldots, n$.  In the context of software testing, Boland et al.\ 
(2002) consider comparisons between $X$ and $Y$ with respect to several stochastic orders.  We note that Theorem 
\ref{discrete} gives an alternative, somewhat faster, derivation of some of their results.  Our starting point is the 
well-known relation $X\leq_{lc} Y$, which is equivalent to Newton's inequalities (Hardy et al.\ 1964).  Thus Theorem 
\ref{discrete} and simple calculations yield 
\begin{enumerate}
\item
$X\leq_{st} Y$ ($X\leq_{hr} Y$) if and only if $p\geq 1-\left(\prod_{i=1}^n (1-p_i)\right)^{1/n}$;
\item
$X\leq_{lr} Y$ ($X\leq_{rh} Y$) if and only if $p\geq 1-n/\left(\sum_{i=1}^n (1-p_i)^{-1}\right)$.
\end{enumerate}
If we let $X'=n-X$ and $Y'=n-Y$, then obviously $X'\leq_{lc} Y'$ and 
\begin{align*}
X'\leq_{st} Y'\Longleftrightarrow Y\leq_{st} X;\quad & X'\leq_{hr} Y'\Longleftrightarrow Y\leq_{rh} X;\\
X'\leq_{lr} Y'\Longleftrightarrow Y\leq_{lr} X;\quad & X'\leq_{rh} Y'\Longleftrightarrow Y\leq_{hr} X.
\end{align*}
Applying Theorem \ref{discrete} to $X'$ and $Y'$, we get 
\begin{enumerate}
\item
$Y\leq_{st} X$ ($Y\leq_{rh} X$) if and only if $p\leq \left(\prod_{i=1}^n p_i\right)^{1/n}$;
\item
$Y\leq_{lr} X$ ($Y\leq_{hr} X$) if and only if $p\leq n/\left(\sum_{i=1}^n p_i^{-1}\right)$.
\end{enumerate}

Our result 
\begin{equation}
\label{correct}
X\leq_{hr} Y\Longleftrightarrow p\geq 1-\left(\prod (1-p_i)\right)^{1/n}
\end{equation}
corrects a slight oversight of Boland et al.\ (2002) (Theorem 1, Part (iv) b).  Basically, Boland et al.\ (2002) find the correct criterion for $Y' \leq_{hr} X'$, and claim that the same criterion holds for $X\leq_{hr} Y$.  However, $Y' \leq_{hr} X'$ is equivalent to $X\leq_{rh} Y$, not $X\leq_{hr} Y$.  This explains the discrepancy between (\ref{correct}) and Theorem 1, Part (iv) b, of Boland et al.\ (2002). 

Theorems \ref{gen} and \ref{discrete} are particularly useful for comparing exponential family distributions with their 
mixtures, as will be illustrated in Section 2, where various specific results concerning Poisson, binomial, negative binomial, 
and gamma distributions are unified and generalized.  Section 3 applies the results of Section 2 to convolutions of gamma 
distributions, which are useful in modeling, for example, the lifetime of a redundant standby system without repairing (Bon 
and Paltanea 1999).  It is shown that, if 
$S=\sum_{i=1}^n \beta_i S_i$, where $S_i\sim {\rm Gam}(\alpha_i, 1)$ independently, $\alpha_i, \beta_i>0$, and 
$T=\beta\sum_{i=1}^n S_i,\ \beta>0$, then 
$$T\leq_{st} S \Longleftrightarrow T\leq_{hr} S \Longleftrightarrow \beta\leq \left(\prod_{i=1}^n  
\beta_i^{\alpha_i}\right)^{1/\alpha_+},$$
where $\alpha_+=\sum_{i=1}^n \alpha_i$.  Moreover, 
$$T\leq_{lr} S\Longleftrightarrow T\leq_{rh} S\Longleftrightarrow  \beta\leq \alpha_+\left/\left(\sum_{i=1}^n \alpha_i/\beta_i\right)\right. .$$
In Section 4 convolutions of negative binomial distributions are considered and results analogous to those of Section 3 are 
obtained. 

\section{Comparing exponential family distributions with their mixtures}
Consider the density of an exponential family
\begin{equation}
\label{exp}
f(x; \theta)=f_0(x)\exp[b(\theta)x]h(\theta),
\end{equation}
where $\theta$ is a parameter, and for simplicity, assume the support of $f(x; \theta)$ is the interval $(0,\infty)$ 
(regardless of the value of $\theta$).  Let $g(x)=\int f(x; t)\, {\rm d}\mu(t)$ be the mixture of $f(x; \theta)$ with respect 
to a probability distribution $\mu$ on $\theta$.  Shaked (1980) considers the comparison between $g(x)$ and 
$f(x; \theta)$ with a fixed $\theta$, focusing on the case when the two distributions have the same mean.  Our comparisons 
here are in terms $\leq_{st},\ \leq_{hr},\ \leq_{rh}$ and $\leq_{lr}$.  As noted by Whitt (1985), 
$$\log(g(x)/f(x; \theta))=\log\left(\int e^{[b(t)-b(\theta)]x}h(t)/h(\theta)\, {\rm d}\mu(t)\right)$$
is a convex function of $x$, i.e., $l(x)=\log(f(x; \theta)/g(x))$ is concave.  
(This holds because log-convexity is closed under mixture.)  We may compute 
$$\lim_{x\downarrow 0} l(x)=-\log \left(\int h(t)/h(\theta)\, {\rm d}\mu(t)\right),$$
and 
$$\lim_{x\downarrow 0} l'(x) =\frac{\int [b(\theta)-b(t)]h(t)\, {\rm d}\mu(t)}{\int h(t)\, 
{\rm d}\mu(t)},$$
provided the interchange of limit (differentiation) and integration is valid.  Thus, if random variables $X$ and $Y$ have 
densities $f(x; \theta)$ and $g(x)$ respectively, then by Theorem \ref{gen}, 
\begin{enumerate}
\item
$X\leq_{st} Y$ ($X\leq_{hr} Y$) if and only if 
\begin{equation}
\label{st}
\int h(t)\, {\rm d}\mu(t)\leq h(\theta);
\end{equation}
\item
$X\leq_{lr} Y$ ($X\leq_{rh} Y$) if and only if
\begin{equation}
\label{lr}
b(\theta)\leq \frac{\int b(t)h(t)\, {\rm d}\mu(t)}{\int h(t)\, {\rm d}\mu(t)}.
\end{equation}
\end{enumerate}

If $f(x; \theta)$ is a discrete pmf on $\mathbf{Z}_+$, then by Theorem \ref{discrete}, 
\begin{enumerate}
\item
$X\leq_{st} Y$ ($X\leq_{hr} Y$) if and only if
\begin{equation} 
\label{std}
\int h(t)\, {\rm d}\mu(t)\leq h(\theta); 
\end{equation}  
\item
$X\leq_{lr} Y$ ($X\leq_{rh} Y$) if and only if
\begin{equation}
\label{lrd}
\exp[b(\theta)]\leq \frac{\int h(t)\exp[b(t)]\, {\rm d}\mu(t)}{\int h(t)\, {\rm d}\mu(t)}.
\end{equation}
\end{enumerate}

Let us illustrate (\ref{std}) and (\ref{lrd}) with some discrete examples.  In Examples 1 and 2, certain results of Misra et 
al.\ (2003) are recovered concerning the comparisons of Poisson and binomial distributions with their mixtures; in Example 3 we 
consider the negative binomial and recover analogous results of Alamatsaz and Abbasi (2008).  In addition to $\leq_{st}$ and 
$\leq_{lr}$ studied by Misra et al.\ (2003) and  Alamatsaz and Abbasi (2008), comparisons in terms of $\leq_{hr}$ and $\leq_{rh}$ are also included.

{\bf Example 1.}  Let $X$ have a Poisson distribution ${\rm Po}(\lambda),\ \lambda>0$, whose pmf is
$$f(x; \lambda)= (1/x!)\lambda^x \exp(-\lambda),\quad x=0,1,\ldots,$$
or, in the form of (\ref{exp}),
$$f(x; \lambda)= (1/x!)\exp(xb(\lambda))h(\lambda),$$
with $b(\lambda)=\log(\lambda)$ and $h(\lambda)=\exp(-\lambda)$.  Suppose $Y$ is a mixture of Po($t$) with respect to a 
distribution $\mu(t)$ on $t\in (0, \infty)$.  Then, by (\ref{std}) and (\ref{lrd}) we have 
\begin{enumerate}
\item
$X\leq_{st} Y$ ($X\leq_{hr} Y$) if and only if
$$\int \exp(-t)\, {\rm d}\mu(t)\leq \exp(-\lambda);$$
\item
$X\leq_{lr} Y$ ($X\leq_{rh} Y$) if and only if
$$\lambda \leq \frac{\int t\exp(-t)\, {\rm d}\mu(t)}{\int \exp(-t)\, {\rm d}\mu(t)}.$$
\end{enumerate}

{\bf Example 2.}  Let $X$ have a binomial distribution with parameters $(n, p)$, where $0<p<1$ and $n$ is a positive integer.  
The pmf of $X$ is 
$$f(x; p)=\binom{n}{x} p^x (1-p)^{n-x},\quad x=0,\ldots, n,$$
or, in the form of (\ref{exp}), 
$$f(x; p)=\binom{n}{x} \exp(xb(p))h(p),$$
with $b(p)=\log(p/(1-p))$ and $h(p)=(1-p)^n$.  Suppose $Y$ is a mixture of binomial($n, t$) with respect to a distribution 
$\mu(t)$ on $t\in (0, 1)$.  Then, after simple algebra, (\ref{std}) and (\ref{lrd}) give
\begin{enumerate}
\item
$X\leq_{st} Y$ ($X\leq_{hr} Y$) if and only if
$$\int (1-t)^n\, {\rm d}\mu(t)\leq (1-p)^n;$$
\item
$X\leq_{lr} Y$ ($X\leq_{rh} Y$) if and only if
$$p\leq \frac{\int t(1-t)^{n-1}\, {\rm d}\mu(t)}{\int (1-t)^{n-1}\, {\rm d}\mu(t)}.$$
\end{enumerate}
By considering $X'=n-X$ and $Y'=n-Y$, we get
\begin{enumerate}
\item
$Y\leq_{st} X$ ($Y\leq_{rh} X$) if and only if
$$\int t^n\, {\rm d}\mu(t)\leq p^n;$$
\item
$Y\leq_{lr} X$ ($Y\leq_{hr} X$) if and only if
$$p\geq \frac{\int t^n\, {\rm d}\mu(t)}{\int t^{n-1}\, {\rm d}\mu(t)}.$$
\end{enumerate}

{\bf Example 3.}  Let $X$ have a negative binomial distribution NB$(k, p)$ where $k$ (not necessarily an integer) is positive 
and $0<p<1$.  The pmf of $X$ is 
$$f(x; p)=\binom{k+x-1}{x} p^k (1-p)^x,\quad x=0,1,\ldots,$$
or, in the form of (\ref{exp}), 
$$f(x; p)=\binom{k+x-1}{x} \exp(xb(p))h(p),$$
with $b(p)=\log(1-p)$ and $h(p)=p^k$.  Suppose $Y$ is a mixture of NB($k, t$) with respect to a distribution $\mu(t)$ on $t\in 
(0, 1)$.  Then (\ref{std}) and (\ref{lrd}) give
\begin{enumerate}
\item
$X\leq_{st} Y$ ($X\leq_{hr} Y$) if and only if
\begin{equation}
\label{nbst}
\int t^k\, {\rm d}\mu(t)\leq p^k;
\end{equation}
\item
$X\leq_{lr} Y$ ($X\leq_{rh} Y$) if and only if
\begin{equation}
\label{nblr}
p\geq \frac{\int t^{k+1}\, {\rm d}\mu(t)}{\int t^k\, {\rm d}\mu(t)}.
\end{equation}
\end{enumerate}

Let us illustrate (\ref{st}) and (\ref{lr}) with a continuous example.

{\bf Example 4.}  Let $X$ have a gamma distribution Gam$(\alpha, \beta),\ \alpha>0,\ \beta>0$, which is parameterized so that 
the pdf is 
$$f(x; \beta)=\Gamma(\alpha)^{-1}\beta^{-\alpha} x^{\alpha-1} \exp(-x/\beta),\quad x>0,$$
or, in the form of (\ref{exp}), 
$$f(x; \beta)=\Gamma(\alpha)^{-1}x^{\alpha-1}\exp(xb(\beta))h(\beta),$$
with $b(\beta)=-\beta^{-1}$ and $h(\beta)=\beta^{-\alpha}$.  Suppose $Y$ is a mixture of Gam($\alpha, t$) with respect to a 
distribution $\mu(t)$ on $t\in (0, \infty)$.  Then (\ref{st}) and (\ref{lr}) give
\begin{enumerate}
\item
$X\leq_{st} Y$ ($X\leq_{hr} Y$) if and only if
\begin{equation}
\label{gamst}
\int t^{-\alpha}\, {\rm d}\mu(t)\leq \beta^{-\alpha};
\end{equation}
\item
$X\leq_{lr} Y$ ($X\leq_{rh} Y$) if and only if
\begin{equation}
\label{gamlr}
\beta\int t^{-\alpha-1}\, {\rm d}\mu(t)\leq \int t^{-\alpha}\, {\rm d}\mu(t)<\infty.
\end{equation}
\end{enumerate}
Note that, unlike previous examples, this is a continuous case and the regularity conditions (interchange of limit 
(differentiation) and integration) required in the derivation of (\ref{gamst}) and (\ref{gamlr}) need to be verified.  For 
example, to establish (\ref{gamst}), we note
\begin{align*}
\lim_{x\downarrow 0} \frac{f(x; \beta)}{\int f(x; t)\, {\rm d}\mu(t)}&=\lim_{x\downarrow 0} \frac{\beta^{-\alpha} 
\exp(-x/\beta)}{\int t^{-\alpha} \exp(-x/t)\, {\rm d}\mu(t)}\\
&=\frac{\beta^{-\alpha}}{\lim_{x\downarrow 0} \int t^{-\alpha} \exp(-x/t)\, {\rm d}\mu(t)}\\
&=\frac{\beta^{-\alpha}}{\int t^{-\alpha}\, {\rm d}\mu(t)},
\end{align*}
where we appeal to the monotone convergence theorem for the last equality.

\section{Convolutions of gamma distributions}
Example 4 in Section 2 enables us to compare a sum of independent gamma random variables with a particular gamma variate.  To 
achieve this we exploit a connection between such a convolution of gamma distributions and a mixture of gamma distributions.
Specifically, let $S=\sum_{i=1}^n \beta_i S_i,$ where $S_i\sim {\rm Gam}(\alpha_i, 1)$ independently and $\beta_i>0,\ i=1, 
\ldots, n$.  Let $T \sim {\rm Gam}(\sum_{i=1}^n \alpha_i, \beta),\ \beta>0$.  We are interested in conditions on $\beta$ that 
ensure $T\leq_{st} S,\ T\leq_{hr} S,\ T\leq_{rh} S$ or $T\leq_{lr} S$.  Relevant works on this problem include Boland et al.\ (1994), Bon 
and Paltanea (1999), Kochar and Ma (1999), Korwar (2002), and Khaledi and Kochar (2004).  In particular, using {\it 
majorization} techniques (Marshall and Olkin 1979), Boland et al.\ (1994) show that, in the case when all $\alpha_i=1$, 
i.e., when $S$ is a sum of independent exponential variables with possibly different scales, we have 
$$\beta\leq \frac{n}{\sum_{i=1}^n \beta_i^{-1}}\ \Longrightarrow\ T\leq_{lr} S.$$
Bon and Paltanea (1999) extend this to (still with $\alpha_i=1$)
\begin{align}
\label{bp1}
T\leq_{st} S\ (T\leq_{hr} S)\ &\Longleftrightarrow\ \beta\leq \left(\prod_{i=1}^n \beta_i\right)^{1/n};\\
\label{bp2}
T\leq_{lr} S\ &\Longleftrightarrow\ \beta\leq \frac{n}{\sum_{i=1}^n \beta_i^{-1}}.
\end{align}
The results of Korwar (2002) and Khaledi and Kochar (2004) imply that the ``$\Longleftarrow$'' parts of (\ref{bp1}) and 
(\ref{bp2}) hold when all $\alpha_i$ are equal, and their common value $\alpha\geq 1$.  As an application of the calculations 
in Sections 1 and 2, we give a further extension for general $\alpha_i>0$.  Such results are of interest in reliability theory 
as they provide convenient bounds (for example) on the hazard rate function of $S$ through the simpler hazard rate function of 
$T$ (Bon and Paltanea 1999).

\begin{theorem}
Assume $\alpha_i>0$ and let $\alpha_+=\sum_{i=1}^n \alpha_i$.  Then
\begin{enumerate}
\item
$T\leq_{st} S$ ($T\leq_{hr} S$) if and only if $\beta\leq \left(\prod_{i=1}^n \beta_i^{\alpha_i}\right)^{1/\alpha_+}$;
\item
$T\leq_{lr} S$ ($T\leq_{rh} S$) if and only if $\beta\leq \alpha_+\left/\left(\sum_{i=1}^n \alpha_i/\beta_i\right)\right.$.
\end{enumerate}
\end{theorem}
\begin{proof} Let $T_0=\sum_{i=1}^n S_i$.  We know that $(S_1/T_0, \ldots, 
S_n/T_0)$ is independent of $T_0$ (property of the gamma distribution); consequently $S/T_0=\sum \beta_i S_i/T_0$ is independent of $T_0$.  Denote the distribution 
of $S/T_0$ by $\mu$.  Then $S=(S/T_0)T_0$ has the distribution of a mixture of ${\rm Gam}(\alpha_+, \gamma)$ with 
respect to $\mu(\gamma)$ on $\gamma\in (0, \infty)$, whereas $T\sim {\rm Gam}(\alpha_+, \beta)$.  Thus the results of 
Example 4, i.e., (\ref{gamst}) and (\ref{gamlr}), are directly applicable.  We only need to calculate 
$$\int \gamma^{-\alpha_+}\, {\rm d}\mu(\gamma)=E[(S/T_0)^{-\alpha_+}]$$ 
and 
$$\int \gamma^{-\alpha_+-1}\, {\rm d}\mu(\gamma) =E[(S/T_0)^{-\alpha_+-1}].$$  

It can be shown that
\begin{align}
\label{diri1}
E[(S/T_0)^{-\alpha_+}]&=\prod_{i=1}^n \beta_i^{-\alpha_i},\ {\rm and}\\
\label{diri2}
E[(S/T_0)^{-\alpha_+-1}]&=\frac{\sum_{i=1}^n \alpha_i/\beta_i}{\alpha_+}\prod_{i=1}^n \beta_i^{-\alpha_i}.
\end{align}
The claims then follow from (\ref{gamst}) and (\ref{gamlr}).  Equation (\ref{diri1}) dates back to Mauldon (1959), and the 
following derivation, which we include for completeness, can be found in Letac et al.\ (2001).  For $t_1, \ldots, t_n\in 
(-\infty, 1)$ we have, by independence, 
\begin{align*}
E\left[\exp\left(\sum t_iS_i\right)\right]&=\prod_{i=1}^n E[\exp(t_iS_i)]\\
&=\prod_{i=1}^n (1-t_i)^{-\alpha_i}.
\end{align*}
On the other hand,
\begin{align*}
E\left[\exp\left(\sum t_iS_i\right)\right]&=E\left\{E\left[\exp\left(\sum t_iS_i\right)|\sum t_iS_i/T_0\right]\right\}\\
&=E\left[\left(1-\sum t_i S_i/T_0\right)^{-\alpha_+}\right].
\end{align*}
Thus 
$$E\left[\left(1-\sum t_iS_i/T_0\right)^{-\alpha_+}\right]=\prod_{i=1}^n (1- t_i)^{-\alpha_i}.$$
Equation (\ref{diri1}) is obtained by substituting $(1-\beta_i)$ for $t_i,\ i=1, \ldots, n$. Moreover, (\ref{diri2}) is 
obtained by differentiating both sides of (\ref{diri1}) with respect to $\beta_i$ and then adding the results for $i=1, 
\ldots, n$. 
\end{proof}

Actually, Khaledi and Kochar (2004) also compare variables of the form of $S$ (assuming $\alpha_i$ are equal and their 
common value $\alpha\geq 1$) in terms of the {\it dispersive order} $\leq_{disp}$.  We mention a result comparing $T$ and $S$ 
in terms of $\leq_{disp}$ for general $\alpha_i>0$.  Let us recall the definitions of $\leq_{disp}$ and the related {\it star 
order} $\leq_*$.
\begin{definition}
Let $X$ and $Y$ be absolutely continuous random variables supported on $(0,\infty)$ with cdfs $F$ and $G$ respectively and 
denote by $F^{-1}$ and $G^{-1}$ the inverse functions of $F$ and $G$ respectively.  
\begin{itemize}
\item
We say $X$ is smaller than $Y$ in the {\it dispersive order}, or $X\leq_{disp} Y$, if 
$$F^{-1}(b)-F^{-1}(a)\leq G^{-1}(b)-G^{-1}(a),\quad 0<a<b<1.$$
\item
We say $X$ is smaller than $Y$ in the {\it star order}, or $X\leq_* Y$, if
$G^{-1} F(x)/x$ is an increasing function of $x,\ x>0$.
\end{itemize}
\end{definition}

\begin{theorem}
\label{disp}
We have $T\leq_{disp} S\Longleftrightarrow T\leq_{st} S$.
\end{theorem}
\begin{proof}  The ``$\Longrightarrow$'' part follows from the definitions (see Theorem 2.B.7 of Shaked and Shanthikumar, 
1994).  To prove the ``$\Longleftarrow$'' part, first we show $T\leq_* S$.  The claim $T\leq_{disp} S$ then follows from 
$T\leq_{st} S$ and $T\leq_* S$ (Ahmed et al.\ 1986; Shaked and Shanthikumar, 1994).  Denote the density functions of $T$ and 
$S$ by $f(x)$ and $g(x)$ respectively.  One sufficient condition for $T\leq_* S$ is that, for all $a>0,\ af(ax)-g(x)$ changes 
sign at most twice as $x$ increases from 0 to $\infty$, the sign sequence being $-,\ +,\ -$ in the case of two changes.  This 
is easily verified by noting that, based on the analysis in Section 2, $\log(af(ax)/g(x))$ is concave in $x$.
\end{proof}
 
\section{Convolutions of negative binomial distributions}
This section contains results for sums of independent negative binomial random variables.  The development somewhat parallels 
that of Section 3. 

Let $N=\sum_{i=1}^n N_i,$ where $N_i\sim {\rm NB}(k_i, p_i)$ independently, $k_i>0,\ p_i\in (0,1),\ i=1, \ldots, n$.  Let 
$M\sim {\rm NB}(\sum_{i=1}^n k_i, p),\ p\in (0,1)$.  For the special case when all $k_i=1$, Boland et al.\ (1994) compare 
variables of the form of $N$, i.e., sums of independent geometric variables with possibly different parameters, with 
respect to the likelihood ratio order.  We have the following result comparing $M$ and $N$ for general $k_i>0$ (not 
necessarily integers).  Theorem \ref{nb} should be compared with Example 0 in Section 1. 

\begin{theorem}
\label{nb}
Let $k_+=\sum_{i=1}^n k_i$.  Then
\begin{enumerate}
\item
$M\leq_{st} N$ ($M\leq_{hr} N$) if and only if $p\geq \left(\prod_{i=1}^n p_i^{k_i}\right)^{1/k_+}$;
\item
$M\leq_{lr} N$ ($M\leq_{rh} N$) if and only if $p\geq \sum_{i=1}^n k_i p_i/k_+$.
\end{enumerate}
\end{theorem}
\begin{proof} The negative binomial NB($k, t$) is a mixture of Po($\lambda (1-t)/t$), where the mixing distribution is 
$\lambda\sim {\rm Gam}(k, 1)$.  It follows that the distribution of $N=\sum_{i=1}^n N_i$ is given by 
\begin{align*}
N|(\lambda_1, \ldots, \lambda_n)&\sim {\rm Po}\left(\sum_{i=1}^n \lambda_i (1-p_i)/p_i\right),\\
                       \lambda_i&\sim {\rm Gam}(k_i, 1)\ {\rm independently}.
\end{align*}
In this setup let $L=\sum_{i=1}^n \lambda_i (1-p_i)/p_i$ and $\lambda_+=\sum_{i=1}^n \lambda_i$.  As in Section 3, 
$L=(L/\lambda_+)\lambda_+$ is a scale mixture Gam($k_+, \gamma$) where the distribution of $\gamma$ is that of $L/\lambda_+$.  
It is clear that $N$ can be expressed as a mixture of negative binomial variates: 
$$N|\gamma \sim {\rm NB}\left(k_+, (1+\gamma)^{-1}\right)$$
where again $\gamma$ has the distribution of $L/\lambda_+$.  We may apply the results of Example 3 in Section 
2, namely (\ref{nbst}) and (\ref{nblr}).  However, as pointed out by an anonymous reviewer, it is simpler to appeal 
to Theorem \ref{discrete} directly.  By the mixture representation of $N$ above we have $M\leq_{lc} N$.   
A quick calculation yields 
$$\frac{\Pr(M=0)}{\Pr(N=0)}=\frac{p^{k_+}}{\prod_{i=1}^n p_i^{k_i}}$$
and
$$\frac{\Pr(M=1)}{\Pr(N=1)}=\frac{k_+ p^{k_+}(1-p)}{\left(\prod_{i=1}^n p_i^{k_i}\right) \sum_{i=1}^n k_i(1-p_i)}.$$
The claims then follow from Theorem \ref{discrete}. 
\end{proof} 

\section*{Acknowledgement}
The author would like to thank the Editor and an anonymous reviewer for their helpful comments.

\end{document}